\documentclass{article}
\usepackage{amstex,amssymb}
\newcommand{\ol}{\setlength{\itemsep}{0pt.}\begin{enumerate}}
\newcommand{\eol}{\end{enumerate}\setlength{\itemsep}{-\parsep}}
\newcommand{\remove}[1]{}

\newtheorem{theorem}{Theorem}

\newcommand{\qed}{\hfill$\Box$}

\newcommand{\ber}{{\begin{eqnarray*}}}
\newcommand{\eer}{{\end{eqnarray*}}}
\newcommand\bysame{\rule[1mm]{1cm}{.025cm}}
\newcommand\nc\newcommand
\nc\bfa{{\bf a}}\nc\bfA{{\bf A}}\nc\cA{{\mathcal A}}
\nc\bfb{{\bf b}}\nc\bfB{{\bf B}}\nc\cB{{\cal B}}
\nc\bfc{{\bf c}}\nc\bfC{{\bf C}}\nc\cC{{\cal C}}
\nc\bfd{{\bf d}}\nc\bfD{{\bf D}}\nc\cD{{\cal D}}
\nc\bfe{{\bf e}}\nc\bfE{{\bf E}}\nc\cE{{\cal E}}
\nc\bff{{\bf f}}\nc\bfF{{\bf F}}\nc\cF{{\cal F}}
\nc\bfg{{\bf g}}\nc\bfG{{\bf G}}\nc\cG{{\cal G}}
\nc\bfh{{\bf h}}\nc\bfH{{\bf H}}\nc\cH{{\cal H}}
\nc\bfi{{\bf i}}\nc\bfI{{\bf I}}\nc\cI{{\cal I}}
\nc\bfj{{\bf j}}\nc\bfJ{{\bf J}}\nc\cJ{{\cal J}}
\nc\bfk{{\bf k}}\nc\bfK{{\bf K}}\nc\cK{{\cal K}}
\nc\bfl{{\bf l}}\nc\bfL{{\bf L}}\nc\cL{{\cal L}}
\nc\bfm{{\bf m}}\nc\bfM{{\bf M}}\nc\cM{{\cal M}}
\nc\bfn{{\bf n}}\nc\bfN{{\bf N}}\nc\cN{{\cal N}}
\nc\bfo{{\bf o}}\nc\bfO{{\bf O}}\nc\cO{{\cal O}}
\nc\bfp{{\bf p}}\nc\bfP{{\bf P}}\nc\cP{{\cal P}}
\nc\bfq{{\bf q}}\nc\bfQ{{\bf Q}}\nc\cQ{{\cal Q}}
\nc\bfr{{\bf r}}\nc\bfR{{\bf R}}\nc\cR{{\cal R}}
\nc\bfs{{\bf s}}\nc\bfS{{\bf S}}\nc\cS{{\cal S}}
\nc\bft{{\bf t}}\nc\bfT{{\bf T}}\nc\cT{{\cal T}}
\nc\bfu{{\bf u}}\nc\bfU{{\bf U}}\nc\cU{{\cal U}}
\nc\bfv{{\bf v}}\nc\bfV{{\bf V}}\nc\cV{{\cal V}}
\nc\bfw{{\bf w}}\nc\bfW{{\bf W}}\nc\cW{{\cal W}}
\nc\bfx{{\bf x}}\nc\bfX{{\bf Z}}\nc\cX{{\cal X}}
\nc\bfy{{\bf y}}\nc\bfY{{\bf Y}}\nc\cY{{\cal Y}}
\nc\bfz{{\bf z}}\nc\bfZ{{\bf Z}}\nc\cZ{{\cal Z}}

\def\dim{\qopname\relax{no}{dim}}
\def\codim{\qopname\relax{no}{codim}}
\def\dist{\qopname\relax{no}{\partial}}

\newcommand\Proof{\noindent{\sc Proof. }}
\newcommand\reals{{\mathbb R}}

\newcommand\ff{{\mathbb F}}
\newcommand\pp{{\mathbb P}}

\begin{document}
\title{Polynomial method in coding and information theory}
\author{A.Ashikhmin
\thanks{Bell Laboratories, Lucent
Technologies, 600 Mountain Avenue, Murray Hill, NJ 07974.
\remove{\tt alexei@c3serve.c3.lanl.gov.}}
\and
A. Barg$^\ast$  
\and
S. Litsyn \thanks{Department of Electrical
Engineering--Systems, Tel Aviv University, Ramat Aviv
69978, Israel, \remove{\tt litsyn@eng.tau.ac.il.}}
}
\date{}
\maketitle\thispagestyle{empty}

\vskip 2cm\begin{abstract}  
Polynomial, or Delsarte's,  method in coding theory
accounts for a variety of structural results on, and bounds
on the size of, extremal configurations (codes and designs) in
various metric spaces. In recent works of the authors the
applicability of the method was extended to cover a wider range
of problems in coding and information theory.
In this paper we present a general framework for the method
which includes previous results as particular cases.
We explain how this generalization leads to new asymptotic bounds
on the performance of codes in binary-input memoryless channels
and the Gaussian channel, which improve the results of Shannon et al.
of 1959-67, and to a number of other results in combinatorial coding theory.
\end{abstract}

\section{Introduction: Some problems of coding and information\\ theory}
Let $X$ be a metric space with distance function
$\dist(\cdot,\cdot)$. A {\em code} $C$ is an arbitrary finite subset of $X$.
The number $d(C)=\min_{\bfc_1,\bfc_2\in C,\,\bfc_1\ne\bfc_2}\dist
(\bfc_1,\bfc_2)$ is called the {\em distance} of $C$. The study of codes was
initiated in the context of transmission of information over
noisy channels \cite{sha48}.

The motivating example is codes in the binary Hamming space
$H^n=\ff_2^n$ with the metric
$\dist(\bfx,\bfy)=|\{e\in\{1,2,\dots,n\}\mid x_e\ne y_e\}|.$
This example corresponds to transmission over the binary
symmetric channel (BSC). Suppose that a vector $\bfx$ is transmitted.
In the channel each coordinate is inverted with probability
$p$ and left intact with probability $1-p$ and different coordinates
are subjected to the error process independently.
Let $P(\bfy|\bfx)$ be the conditional probability distribution
induced on $H^n$ by this channel. $P(\bfy|\bfx)$
is a monotone (decreasing) function of the distance $\dist(\bfx,\bfy)$;
hence it is possible to study the performance of codes in geometric
terms. In particular, for small $p$ the most important parameter
of the code is its distance. This gives an information-theoretic reason
to look for codes of a given size with large distance.
There are also other combinatorial and geometric reasons for
this interest; we outline them below.

Now suppose that $\bfx\in\reals^n$ and the error process in the channel is
described as follows: for a transmitted vector $\bfx$ we
receive from the channel a vector $\bfy=\bfx+\bfe$, where
each coordinate of $\bfe$ is a Gaussian $(0,\sigma^2)$ random variable.
A consistent definition of capacity of such a channel is obtained
if the input signals satisfy some sort of energy
constraints. Typically one assumes that the energy, or the average energy of
input signals does not exceed a given number $A\sigma^2$ per
dimension, where
$A$ is a positive number called the ``signal-to-noise ratio.''
Shannon \cite{sha59} has shown that for a set of input
signals of sufficiently large size the study of the channel
is reduced to considering signals of {\em constant} energy equal to
$\sigma \sqrt{An}$, that is, points on the $n$-dimensional
sphere. Thus our second example will be $X=S^{n-1}(\reals)$, the
unit sphere in $\reals^n,$ with the Euclidean distance
$\dist(\bfx,\bfy)=||\bfx-\bfy||^{1/2}.$

The third standard example is $X=\{\bfx\in H^n| \#\{i:x_i=1\}=v\}$,
called the binary Johnson space $J^{n,v}$.
Since this space is a subset of $H^n$, it is again associated
with transmission over the BSC.
However, the interest in the Johnson space is largely determined
by the fact that combinatorially it can be studied by methods
similar to the Hamming case \cite{del73}, \cite{ban84} and has
strong connections to the latter \cite{mce77a}, \cite{rod80}.

The theory and a part of results outlined in this paper are valid in a
large class of finite spaces that afford the structure of an
association scheme, and in the infinite case, in compact two-point homogeneous
spaces. However below we concentrate on the above examples since they
give rise to central asymptotic problems of coding and 
information theory that can be treated in geometric terms. Let us 
outline these problems. We use the mixed entropy and entropy functions
\begin{align*}
T_s(x,y)&=x\log_q(q-1)-x\log_s y-(1-x)\log_s(1-y),\\
H_s(x)&=T_s(x,x),
\end{align*}
and omit the subscript if the logarithms are taken base $e$.

\subsection{The size-distance ($R$-$\delta$) problem.}
Let $A(X;n,d)=\max_{C\in X}\{|C|\mid d(C)=d\}$ be the maximal
possible size of the code with a given distance. Finding this function
is one of the central problems of coding theory. Apart from a number
of particular cases for small $n$ this problem is unsolved.

First let $X=H^n$. Let $R=(1/n)\log_2 |C|$ and $\delta=d(C)/n$ 
be the {\em rate} and the {\em relative distance} of the code.
Clearly, $0\le R\le 1, 0\le\delta\le 1.$ Let
\[
\bar R(\delta)=\limsup_{n\to\infty}\frac{\log_2 A(H^n;n,d)}{n},
\ \ \ \underline R(\delta)=\liminf_{n\to\infty}\frac{\log_2 A(H^n;n,d)}{n},
\]
where the limits are computed over all sequences of codes $C_n$
for which $\limsup\limits_{n\to\infty} \frac{d(C_n)}n\ge \delta.$
Below we assume that these two functions have a common limit, denoted
$R(\delta)$ (if they do not, the upper bounds become bounds on
$\bar R$ and the lower ones on $\underline R$).
It is clear that $R(\delta)$ is a monotone decreasing 
function of $\delta;$ its inverse is denoted below by $\delta(R).$

It is known and easily proved that $R(\delta)=0$ for $\delta\in[\frac 12,1].$ 
Otherwise the best known bounds on $\delta(R)$ have the form:
\begin{align}
\delta(R)&\ge \delta^{\rm (vg)}(R):=H_2^{-1}(1-R)\quad
{\mbox{\cite{gil52}, \cite{var57}}}
\label{eq:vg}\\
\delta(R)&\le \delta^{\rm (lp)}(R):=
\min_{\begin{Sb}
0 \le \beta \le \alpha\le 1/2\\ 
H_2(\alpha)-H_2(\beta)=1-R
\end{Sb}}
2\frac{\alpha(1-\alpha)-\beta (1-\beta )} 
{1+2\sqrt{\beta (1-\beta)}}. \quad{\mbox{\cite{mce77a}}}\label{eq:mrrw}
\end{align}

Let $X=S^{n-1}(\reals).$ In this case the corresponding functions are
written in the form $R(d), d(R),$ where $d(R), 0\le d(R)\le2,$ is 
the limit value of the Euclidean distance of codes of rate 
$R:=\frac 1n\log |C|.$ We have 
\begin{align}
d(R)&\ge d^{\rm (s)}(R):=\sqrt{2(1-\sqrt{1-e^{-2R}})} \quad(0\le R<\infty)
\quad{\mbox{\cite{sha59}}} \label{eq:vg-sphere}\\
d(R)&\le d^{\rm (kl)}(R):=\frac{\sqrt2(\sqrt{1+\rho}-\sqrt{\rho})}
{\sqrt{1+2\rho}}
\quad{\mbox{\cite{kab78}}},\nonumber
\end{align}
where in the last formula $\rho$ is the root of 
\begin{equation}\label{eq:rho}
R=(1+\rho)H\big(\frac{\rho}{1+\rho}\Big) \quad(0\le R<\infty).
\end{equation}

Lower bounds (\ref{eq:vg}) and (\ref{eq:vg-sphere}) are obtained by 
random choice. The upper bounds are obtained by the polynomial method 
which is outlined in the next section.

\subsection{Error probability of decoding.}
Though the packing $(R$-$\delta)$ problem has received more attention,
Shannon's original motivation was reliable transmission of information
over channels. Let $X=H^n$ and let $C\subset X$ be a code
used for transmission over the binary symmetric channel.
A {\em decoding} is a (partial) mapping $\psi:X \to C.$ Let $S(t,\bfc)$
be a sphere of radius $t$ around a point $\bfc.$
Consider the decoding defined on $\cup_{\bfc\in C}S(t,\bfc)$ as follows:
\[
\psi_t(\bfy)=\bfc \quad\mbox{if } \bfy\in S(t,\bfc)
\mbox{ and } \dist(\bfy,\bfc)\le\dist(\bfy,\bfc') \mbox{ for all }\bfc'\in C.
\]
Clearly, if $t\le\lfloor(d(C)-1)/2\rfloor$, the decoding result is defined
uniquely; otherwise
we agree that ties are broken arbitrarily. As long as
$t\le\lfloor(d(C)-1)/2\rfloor$,
the spheres $S(t,\bfc)$ are disjoint; for larger $t$ some of them intersect.
Starting with a certain value of $t$ (called the {\em covering radius} of $C$)
their union covers the entire $X$. In this case the decoding
is called {\em complete}. Let $P_{de}(H^n;C,p)$ be the average
error probability of complete decoding for a code $C$ used over
the BSC with error probability $p$:
\[
P_{de}(H^n;C,p):=\frac{1}{|C|}\sum_{\bfx\in C}P_{de}(\bfx),
\]
where the last probability describes the event that the transmitted
code vector is $\bfx$ and the decoding result is a code vector $\bfx'\ne\bfx.$
Let
\[
P_{de}(H^n;n,R,p)=\max P_{de}(H^n;C,p),
\]
where the maximum is taken over all codes of rate $\ge R.$

Obviously, these definitions are valid for any metric space; in
particular for $S^{n-1}(\reals)$ and the Gaussian channel. Therefore,
we also consider the error probability of decoding
$P_{de}(S^{n-1};n,R,A),$ defined analogously.
It is known that $P_{de}(X;n,R)$ falls exponentially for
both $X=H^n$ \cite{eli55b} and $X=S^{n-1}$ \cite{sha59}; consider 
therefore the exponents
\begin{eqnarray*}
\bar E_{de}(H^n;R,p)&=&\limsup\limits_{n\to\infty}-
\frac 1n\log_2 P_{de}(H^n;R,p,n)\\
\bar E_{de}(S^{n-1};R,A)&=&\limsup\limits_{n\to\infty}-\frac 1n
\log P_{de}(S^{n-1};R,A,n).
\end{eqnarray*}
After Shannon \cite{sha59} the best attainable error exponent is
called the {\em reliability function} of the channel. 
Computing the reliability function of these and other channels
dominated information theory through the end of the 1960s \cite{gal68}.
Even in the simplest cases mentioned this problem is still unsolved.
Upper bounds on $\bar E_{de}(H^n;R,p)$ were derived in \cite{sha67}; see also
\cite{mce77b}. Lower bounds on $\underline E_{de}(H^n;R,p)$ were given in 
\cite{eli55b}, \cite{gal65}.
\footnote{Lower bounds on the reliability function constitute Shannon's
channel coding theorem, proved in the general case by Feinstein \cite{fei54}, 
see also Khinchin \cite{khi57}.}
Lower and upper bounds on $E_{de}(S^{n-1}(\reals);R,p)$
were obtained in \cite{sha59}; see also \cite{kab78}.

For $X=H^n$ coding theorists have also studied the
other limiting case of decoding, that of decoding radius $t=0$,
called {\em error detection}.
The probability of undetected error is defined analogously:
\[
P_{ue}(H^n;C,p):=\frac{1}{|C|}\sum_{\bfx\in C}P_{ue}(\bfx),
\]
where the probability $P_{ue}(\bfx)$ corresponds to the event
that the received vector, equal to
$\bfx+\bfe, \bfe\in H^n\!\setminus\!{\bf0},$
is itself in $C$.

\section{Polynomial method}

Delsarte \cite{del72b} suggested a method of deriving
upper bounds on the size of a code with a given distance by optimizing
a certain functional on the cone of polynomials of degree at
most $n$. The formalism of the method can be developed either
in the context of association schemes \cite{del73}, \cite{ban84} or
of harmonic
analysis on noncommutative compact groups \cite{kab78}, \cite{lev83a}.

We again begin with the binary Hamming space $H^n$.
The main role in the method is played by the distance distribution
of codes. Let $C\subset H^n$ be a code. Its (average) {\em distance
distribution} is given by $\cA=(\cA_i, 0\le i\le n)$, where
\[
\cA_i=\frac{1}{|C|}|\{(\bfc,\bfc')\in C^2\mid \dist(\bfc,\bfc')=i\}|.
\]

Let $X=H^n$ and let $(K_k(x), k=0,1,\dots)$ be the system of Krawtchouk
polynomials \cite{sze75}, i.e., polynomials orthogonal on the
set $(0,1,2,\dots,n)$ with weight $\mu(i)=\binom{n}{i}2^{-i}.$
Let $C$ be a code with the distance distribution $\cA.$ The
{\em MacWilliams transform} of $\cA$ is a vector
$\cA'=\frac 1{|C|}\cA\bfK,$ where
$\bfK=(K_k(i), 0\le i,k\le n)$ is the $(n+1)\times(n+1)$ Krawtchouk
matrix with rows numbered by $i$ and columns by $k$. 
One of the central properties of
codes, the {\em Delsarte
inequalities} \cite{del72b},\,\cite{del73}, is that the components of
$\cA'$ are nonnegative.

Now let $f(x)=\sum_{i=0}^n f_iK_i(x)$ be a polynomial of degree at most $n$
written in the Krawtchouk basis. The definition of the MacWilliams
transform implies the following useful identity \cite{del72b}:
\begin{equation}\label{eq:Fourier}
|C| \sum_{i=0}^n f_i\cA'_i=\sum_{i=0}^n f(i)\cA_i.
\end{equation}

The following theorem, proved in particular cases in \cite{del72b},
\cite{ash99a}, \cite{lit99}, \cite{ash99f}, is the main general result
of this paper. It accounts for the new upper bounds on the reliability
functions of the next sections as well as for some other estimates
of code parameters.
\begin{theorem} \label{thm:main}
Let $C$ be a code with distance distribution $\cA.$
Let $f(x)=\sum_{k=0}^n f_kK_k(x),$ $f_k\ge 0, 1\le k\le n,$
be a polynomial of degree at most $n.$ Let $F=\sum_{i=1}^n g(i)\cA_i$
be a function on $C$ and suppose that $f(i)\le g(i), 0\le i\le n.$
Then
\[
F\ge |C|f_0-f(0).
\]
\end{theorem}
\Proof By the inequality $\cA'\ge 0$ and (\ref{eq:Fourier}) we obtain
\begin{equation*}
|C|f_0\le |C|\sum_{i=0}^n f_i\cA'_i=f(0)+\sum_{j=1}^n f(j)\cA_j
\le f(0)+\sum_{j=1}^n g(j)\cA_j=f(0)+F.
\end{equation*}
where we have used the fact that $\cA_0'=1.$
\qed

{\em Examples}.

1. {Probability of undetected error.}
Let $0\le p\le \frac 12$ and $g(i)=p^i(1-p)^{n-i}.$ Then
$F=P_{ue}(H^n;C,p).$

2. {Delsarte's linear programming bound.}
Let $g(i)=0, 1\le i\le n$, then $F=0$.
Suppose that the code $C$ in Theorem 1 has distance $d$. Then
$\cA_i=0$ for $i=1,\dots,d-1.$ Hence it suffices to assume that
$f(i)\le 0$ for $i=d,d+1,\dots,n.$ Assuming in addition that
$f_0>0,$ we obtain Delsarte's linear programming bound on the size of a
code with distance $d$ \cite{del72b}:
\begin{equation}\label{eq:lp-bound}
|C|\le\inf_f\big\{\frac{f(0)}{f_0}\big|
f_0>0, f_k\ge 0, 1\le k\le n; f(i)\le 0, i=d,d+1,\dots,n\big\}.
\end{equation}
The problem of finding stationary points of the functional $f(0)/f_0$
has been one of the central in combinatorial coding theory since 1972
(see \cite{lev98}).

3. Let $1\le w\le n$ be an integer and let $g(i)=\binom{n-i}{n-w}$.
We obtain a set of code invariants
$F_w=\sum_{i=0}^w\binom{n-i}{n-w} \cA_i$.
The numbers $F_w$ (binomial moments of the distance distribution)
are related to numerous combinatorial invariants
\cite{ash99a}, \cite{bar99a}, \cite{bar99b}, for instance, the cumulative size
of subcodes of restricted support, and, in the linear case, to the
higher weight enumerators, rank polynomial, Tutte polynomial, etc.

4. Suppose that in Theorem \ref{thm:main} $f(i)> 0$ for $0\le i\le w$
and $f(i)\le 0$ for $w+1\le i\le n$, where $w\in[1,n]$ is a
parameter. Put $g(i)=f(i).$
Then the theorem implies in a way similar to Example 2 the inequality
\[
\sum_{i=1}^w f(i)\cA_i\ge |C|f_0-f(0).
\]
In other words, there exists a number $j,1\le j\le w,$ such that
\begin{equation}\label{eq:lower-spec}
\cA_j\ge\frac{|C|f_0-f(0)}{f(j)}.
\end{equation}
This is one of the main results in \cite{lit99}. Since the polynomial
$f$ has to satisfy the same conditions as in Example 2, it is possible
to use the known results in the $R$-$\delta$ problem to derive
specific lower bounds on the distance distribution of codes \cite{lit99}
(see Sect.\ref{sec:lower}).

Let $C\subset J^{n,v}$ be a code in the Johnson space
and $\cA=(\cA_{2i},0\le i \le v)$ be its average distance distribution.
Let $Q^v_k(x)$ be a family of Hahn polynomials orthogonal on
the set $(0,1,\dots,v)$ with weight $\mu(i)=\tfrac{\binom
vi\binom{n-v}i}{\binom nv}.$ Then as above one can consider
the transformed distribution $\cA'=\frac 1{|C|}\cA{\bf Q},$ where
${\bf Q}=(Q^v_k(i), 0\le i,k\le v)$ is the Hahn matrix. Again by
Delsarte's theory \cite{del73} the components of $\cA'$ are nonnegative.
Thus, one can consider the invariants of Examples 1-4 for the Johnson space.
In particular, inequalities (\ref{eq:lp-bound}),(\ref{eq:lower-spec})
are straightforward \cite{del73},\cite{lit99}, where this time $f(x)$
is a polynomial with positive Hahn-Fourier coefficients.

One particular reason to study bounds of the form (\ref{eq:lp-bound}),
(\ref{eq:lower-spec}) in the Johnson space follows from the fact that
one can translate them to the Hamming space. Namely, since
$J^{n,v}\subset H^n$, an easy averaging argument (the {\em multiple packing
principle} \cite{bas65}) shows that the maximum sizes of codes in $H^n$
and $J^{n,v}$ are related as follows:
\begin{equation}\label{eq:be}
\binom nv A(H^n;n,d)\le 2^n A(J^{n,v};n,d).
\end{equation}
Thus, upper bounds on $A(J^{n,v};n,d)$ also give upper bounds
on codes in the Hamming space (an important example being (\ref{eq:mrrw})).
Several generalizations of this argument are known; in particular,
one can prove a lower bound of the form (\ref{eq:lower-spec}) in
$J^{n,v}$ and then translate it to $H^n$ \cite{lit99}.
An alternative approach to (\ref{eq:be}) is based on
the fact that positive definite functions in the Hahn basis
are also positive definite in the Krawtchouk basis \cite{rod80}. This gives
an analytic method of deriving inequalities of the type (\ref{eq:be}),
which is useful
for those code invariants that are not well defined in the Johnson
space, hence do not carry a geometric meaning. For instance,
this is the case with the $F_w$-invariants of Example 3 \cite{ash99a}.

Finally, consider the case $C\subset S^{n-1}(\reals).$ Let
$\dist(\bfc,\bfc')=||\bfc-\bfc'||$ be the Euclidean distance in
$\reals^n$. It is convenient to define the distance distribution of
$C$ with the help of the function $t(x)=1-\frac{x^2}2.$ In particular,
$t(\dist(\bfc,\bfc'))=\langle\bfc,\bfc'\rangle,$ where
$\langle\cdot,\cdot \rangle$ is the scalar product in $\reals^n$.
Let
\[
a(s,t):=\frac1{|C|}|\{(\bfc,\bfc')\in C^2: s\le \langle\bfc,\bfc'\rangle\le t
\}|
\]
be the {\em distance density} of $C$. Delsarte's inequalities in
this case take on the form \cite{kab78}
\[
\int_{-1}^1\delta(t)a(t,t)P_k^{\lambda,\lambda}(t)dt=\sum_{\bfc,\bfc'\in C}
P_k^{\lambda,\lambda}(\langle \bfc,\bfc'\rangle)\ge 0, \quad k=0,1,\dots,
\]
where $P_k^{\alpha,\beta}(x)$ is the Jacobi polynomial, $\delta(t)$
is the delta-function, and $\lambda=(n-3)/2.$
The analog of (\ref{eq:lower-spec}) in this case is given by the
following theorem.
\begin{theorem}\label{thm:main-sphere}{\rm \cite{ash99f}} Let $C\subset
S^{n-1}(\reals)$
be a code and let $m$ be a integer. Let $-1\le u_0<t(d(C))$ and
suppose that $u_0<u_1<\dots<u_{m-1}<u_m=t(d(C))<1$ are the defining
points of a partition of the segment $[u_0, t(d(C))]$ into $m$
equal segments $U_i=[u_i,u_{i+1}].$

Suppose that $f(x)=\sum_{k=0}^l
f_kP_k^{\alpha,\alpha}(x)$ is a polynomial of degree $l$ such that
$f_k\ge 0, 1\le k\le l,$ and $f(x)\le 0, -1\le x\le u_0,$
$f(x)\ge 0, u_0\le x\le 1.$ Then
there exists a number $i,0\le i\le m-1,$ and a point $s\in U_i$
such that
\begin{equation}\label{eq:lower-spectrum-spherical}
a(u_i,u_{i+1})\ge\frac {f_0|C|-f(1)}{mf(s)}.
\end{equation}
\end{theorem}

The three metric spaces considered above (and many other spaces)
can be studied from one and the same point of view.
This is the principal achievement of \cite{kab78}.
It turns out that the polynomials associated with the space
($K_k, H_k, P_k^{\alpha,\beta}$) represent the zonal spherical kernels
that arise in the analysis of irreducible unitary representations
of the isometry group of the space. Spaces in which
zonal spherical functions are expressed by univariate polynomials
are sometimes called {\em polynomial} \cite{lev83a}, \cite{god93}.

\section{Asymptotics of orthogonal polynomials}

To derive asymptotic bounds on the distance distribution of codes
and other invariants we need asymptotic formulas for orthogonal
polynomials involved in inequalities (\ref{eq:lower-spec}),
(\ref{eq:lower-spectrum-spherical}). These problems have been
studied more or less independently in coding theory
\cite{mce77a}, \cite{kab78}, \cite{lev83a}, \cite{kal95a}, \cite{lit99},
\cite{ash99a}, \cite{ash99f} and analysis \cite{moa79}, \cite{che91},
\cite{gaw91}, \cite{elb94}, \cite{dra97}, \cite{kui99a}.
We quote results from the coding-theory side since they are in the
form better suited to our needs.

Asymptotics of extremal zeros found in \cite{mce77a},\cite{kab78} were
used in these papers to
derive the bounds $\delta^{\rm (lp)}(R)$ and $d^{\rm (kl)}(R)$,
respectively. However, to derive bounds on code invariants we
need to find the behavior of the polynomials from the extremal zero
to the end of the orthogonality segment.

{\em Krawtchouk polynomials}.
$K_k(0)=\binom nk$ and the polynomial is monotone decreasing
in the segment $[0,x_1(K_k)]$, where $x_1$ is the smallest zero of $K_k(x).$
Let $k/n\to\tau$ as $n\to\infty.$
It is known \cite{mce77a}, \cite{lev83a} that
$x_1(K_k)\approx \frac n2-\sqrt{k(1-k)}.$
An asymptotic expression for the exponent of $K_k(x)$ for $x\in[0,x_1(K_k)]$
was derived in \cite{kal95a}. It has the form
\begin{equation}\label{eq:est-Krawtchouk-exact}
\frac 1n\log_2 K_k(\xi n)=H_2(\tau)+
\int_0^\xi \log_2\frac{1-2\tau+\sqrt{(1-2\tau)^2-4y(1-y)}}{2-2y}dy+o(1).
\end{equation}

{\em Hahn polynomials}. The smallest zero of $Q^v_k(x)$ behaves as follows
\cite{mce77a}, \cite{lev83a}:
\[
x_1(Q^v_k)\approx \frac{v(n-v)-k(n-k)}{n+2\sqrt{k(n-k)}}.
\]
Similarly to \cite{kal95a} we have \cite{ash99a},\cite{lit99}
\begin{multline}\label{eq:est-Hahn-exact}
\frac 1n\log_2 Q_k^v(\xi n)=
H_2(\beta)+\int_0^\xi\log_2\Big[\frac{\alpha(1-\alpha)-y(1-2y)-\beta(1-\beta)}
{2(\alpha-y)(1-\alpha-y)}\\
+\frac{\sqrt{[\alpha(1-\alpha)-y(1-2y)-\beta(1-\beta)]^2-4(\alpha-y)(1-\alpha-y)y^2}}
{2(\alpha-y)(1-\alpha-y)}\Big]dy+o(1),
\end{multline}
where $n\to\infty, v=\alpha n, k=\beta n, x\in[0,x_1(Q_k^v)].$

{\em Jacobi polynomials}. The asymptotic expression for the
largest zero of $P_k^{ak,bk}$ has the form \cite{kab78}, \cite{moa79}
\[
x_1^{a,b}:=x_1(P_k^{ak,bk})
\approx \frac{4\sqrt{(a+b+1)(a+1)(b+1)}-a^2-b^2}{(a+b+2)^2}.
\]
The smallest zero then is $-x_1^{b,a}.$
The asymptotic behavior of the exponent of $P_k^{ak,bk}, k\to\infty,$ 
in the entire orthogonality segment
was found in \cite{ash99f}. We quote one of the results in \cite{ash99f}:
let $x\in [-1,-x_1^{b,a}-\epsilon_k]\cup[1,x_1^{a,b}+\epsilon_k],$ 
where $\epsilon_k=k^{-\gamma}, 0\le\gamma\le 1/2.$ Then
\begin{multline}\label{eq:Kalai}
\frac 1k\ln |P_k^{\alpha,\beta}(x)|=
(1+a)H\big(\frac{a}{1+a}\big)\\
\mp\int_x^1 \frac{(a+(a+b)z-b)\mp\sqrt{(a+(a+b)z-b)^2-4(1-z^2)(1+a+b)}}
{2(1-z^2)}dz+o(1),
\end{multline}
where the $-$ sign corresponds to the left of the 2 segments in the
domain of $x$ and the
$+$ to the right of them.

\section{Lower bounds on code invariants}\label{sec:lower}

Specifications of Theorem \ref{thm:main} enable one to prove a large
number of results on code properties. In this section we present
estimates on the distance distribution of codes and related invariants.

Let $C\subset H^n$ be a code of rate $R$ and $\cA$ its distance distribution.
Theorem \ref{thm:main} together with (\ref{eq:est-Krawtchouk-exact})
implies the following
\begin{theorem} \label{thm:spectrum-Hamming} {\rm\cite{lit99}}
For any code of sufficiently large length $n$
there exists a number
$\xi\in[0,1/2-\sqrt{\tau (1-\tau)}], 0\le \tau\le H_2^{-1}(R),$
such that
\[
\cA_{\lfloor \xi n\rfloor}\ge R-H_2(\tau)-2I(\xi,\tau)-o(1),
\]
where $I(\xi,\tau)$ is the integral on the right-hand side of
{\rm(\ref{eq:est-Krawtchouk-exact})}.
\end{theorem}
The bound in this theorem can be slightly improved with the
help of a generalization of (\ref{eq:be}) and asymptotics
(\ref{eq:est-Hahn-exact}) \cite{lit99}.

The lower estimate on $F_w$ invariants of $C$, which is
also proved with the help of Theorem \ref{thm:main}, has the
following form.
\begin{theorem}\label{thm:binom} {\rm\cite{ash99a}}
For any code of sufficiently large length $n$
\[
\frac 1n \log_2 F_{\lfloor 2\omega n\rfloor}\ge
R-1+H_2(\omega^\ast)+(1-\omega^\ast)H_2\Big(\frac{1-2\omega}{1-\omega^\ast}
\Big)-o(1),
\]
where
\[
\omega^\ast=\begin{cases} \omega, &\delta^{{\rm(lp)}}(R)\le\omega\le 1\\
\delta^{{\rm(lp)}}(R), &\delta^{{\rm(lp)}}(R)/2\le\omega\le
\delta^{{\rm(lp)}}(R),
\end{cases}
\]
and $\delta^{{\rm(lp)}}(\cdot)$ is defined in {\rm(\ref{eq:mrrw}).}
\end{theorem}

For $X=S^{n-1}(\reals)$ Theorem \ref{thm:main-sphere} and (\ref{eq:Kalai})
imply the following
\begin{theorem} \label{thm:spectrum-sphere} {\rm\cite{ash99f}}
Let $C$ be a code of rate $R$. Let $\gamma \in[0,\rho],$
where $\rho$ is the root of {\rm (\ref{eq:rho})} be a fixed number.
Then there exists a number $x, \dfrac{2\sqrt{\gamma(1+\gamma)}}{1+2\gamma}
\le x\le 1,$ such that for sufficiently large $n$
\begin{align*}
\frac 1n \ln a(x,x+\frac 1n)\ge 4\gamma(1+\gamma)\int_x^1
\frac {dz}{z+\sqrt{z^2-4(1-z^2)\gamma(1+\gamma)}}-(1+\gamma)H\Big(
\frac \gamma{1+\gamma}\Big)+R-o(1).
\end{align*}
\end{theorem}

These estimates lead to new upper bounds on the reliability function
of the BSC \cite{lit99}, the Gaussian channel \cite{ash99f},
of the exponent of error detection \cite{ash99a}, \cite{ash99f}
(see the next section) and on a number of other parameters of codes.
The approach developed in \cite{kab78} enables us to derive
similar bounds on the distance distribution of codes in projective
real and complex spaces \cite{ash99f}.

\section{Reliability functions and error detection}

The key to the results of this section is given by the following
observation: if a code vector that has many close neighbors is sent over
the channel, the error probability of decoding cannot be too low.
Together with the estimates of the previous section and some
other combinatorial and geometric considerations this leads
to the following results.
\begin{theorem}\label{thm:rel-BSC}{\rm\cite{lit99}}
The reliability function of BSC with error probability $p$ satisfies the 
upper bound 
\[
\bar E_{de}(H^n;R,p)  \le  \max_{\alpha,\beta,\xi,\delta} E_{\alpha,
\beta,
\xi,\delta},
\]
where
\[
E_{\alpha,\beta,\xi,\delta}=\min \Big(
-\delta \log_2 \sqrt{4p(1-p)},
-\tilde{\nu}-\xi \log_2 \sqrt{4p(1-p)} \Big),
\]
and $\alpha$, $\beta$, $\delta,$ and $\xi$ are such that
\(
0 \le \beta \le \alpha \le 1/2,\quad
H_2(\alpha)-H_2(\beta) = 1-R, \ \delta \in [0,\delta^{\mbox(lp)}(R)],
\)
$$\xi \in
\left[ 0,2 \frac{\alpha(1-\alpha)+\beta(1-\beta)}{1+2 \sqrt{\beta
(1-\beta)}}
 \right]; 
$$
\begin{multline*}
\tilde{\nu}=
\min \Big( \nu,
\xi+(1-\xi) H_2(p)
-\max_{\eta \in [ \delta p/2, \min (\delta/4, p(1-\xi))]}
\Big( \delta H_2 \Big(  \frac{2 \eta}{\delta} \Big) 
\\
+(\xi-\delta/2) H_2 \Big(  \frac{\xi-2\eta}{2 \xi -\delta} \Big)
+
(1-\xi-\delta/2)
H_2 \Big( \frac{p(1-\xi)-\eta}{1-\xi-\delta/2} \Big)
\Big),
\end{multline*}
$$\nu=R-1 + H_2(\beta)+2 H_2(\alpha)-
2 q(\alpha,\beta,\xi/2)-\xi-(1-\xi)
H_2\Big(\frac{\alpha-\xi/2}{1-\xi}\Big), $$
and $q(\alpha,\beta,\xi)$ is the function on the right-hand side of
{\rm (\ref{eq:est-Hahn-exact})}.
\end{theorem}

As shown in \cite{sha67}, given any convex upper bound on $E(H^n;R,p)$
one can draw a common tangent to it and the {\em sphere-packing bound}
(one of the bounds in \cite{sha67}), and the segment between the
tangency points will also give an upper bound on $E(H^n;R,p)$, the
so-called
{\em straight-line bound}. Together with the last theorem
this gives the best upper bound to-date on $E(H^n; R,p)$.
Standard methods of information theory \cite{gal68} enable one
to extend this result to memoryless channels with binary input
alphabet.

Error-correcting properties of codes on $S^{n-1}(\reals)$
are given by the following theorem, whose proof relies, in
particular, on Theorem \ref{thm:main-sphere} and (\ref{eq:Kalai}).
\begin{theorem} {\rm \cite{ash99f} }
The reliability function of the
Gaussian channel with signal-to-noise ratio $A$ satisfies the upper
bound
\begin{equation}\label{eq:new-bound}
\bar E_{de}(S^{n-1};R,A)\le\min_{0\le\gamma\le \rho}
\max_{w,d}
\Big[\min (A\frac{d^2}8,A\frac{w^2}8-\cL(w,d,\gamma))\Big],
\end{equation}
where
\[
0\le d\le \frac{\sqrt2\,\big(\sqrt{1+\rho}-\sqrt\rho\big)}{\sqrt{1+2\rho}},
\quad  d\le w\le \frac {\sqrt2\,\big(\sqrt{1+\gamma}-
\sqrt\gamma\big)}{\sqrt{1+2\gamma}},
\]
$\rho$ is the root of {\rm (\ref{eq:rho})}
\begin{gather*}
\cL(w,d,\gamma)=\min\Big\{\frac{Ad^2w^2}{8(4w^2-d^2)},
F\big(1-\frac 12 w^2,\gamma\big)\Big\},\\
F(x,\gamma)=R-(1+\gamma)H\Big(\frac{\gamma}{1+\gamma}\Big)+\int_{x}^1
\frac{4\gamma(1+\gamma)dz}{z+\sqrt{z^2-4(1-z^2)\gamma(1+\gamma)}}.
\end{gather*}
\end{theorem}

The remark made after Theorem \ref{thm:rel-BSC} regarding the
straight-line bound is valid for the Gaussian channel as well; 
taken together, this, again, is the best result known to-date.

Results of the previous sections also lead to the following upper
bound on the exponent of error detection $\bar E_{ue}(H^n;R,p)$ for codes on
the BSC with crossover probability $p$.
\begin{theorem}
\begin{equation*}\label{eq:Pue-LP2}
\bar E_{ue}(H^n;R,p)\le\begin{cases}
1-R-H_2(\delta^{\rm (lp)}(R))+T_2(\delta^{\rm (lp)}(R),p),
&0\le R\le R^{\rm (lp)}(p)
\\
1-R,& R^{\rm (lp)}(p)\le R\le 1,
\end{cases}
\end{equation*}
where $\delta^{\rm (lp)}(\cdot)$ is given by {\rm(\ref{eq:mrrw})}
and $R^{\rm (lp)}(\cdot)$ is its inverse function.
\end{theorem}
This theorem was proved in \cite{ash99a} via lower bounds on the
$F_w$-invariants (Theorem \ref{thm:binom}) and in \cite{lit99} with the use of
Theorem \ref{thm:spectrum-Hamming}. Together with known lower bounds
(see, e.g., \cite{lev89}), it shows that the function $E_{ue}(H^n;R,p)$ is
known exactly for $R\in[ R^{\rm (lp)}(p),1]$.

It is worth mentioning that if the Varshamov-Gilbert bound
(\ref{eq:vg}) is tight (as is widely believed) then
the known lower bounds on $E_{de}(H^n;R,p)$ and $E_{ue}(H^n;R,p)$ are
also tight. The same is true with respect to the Shannon bound
(\ref{eq:vg-sphere}) and $E_{de}(S^{n-1};R,A)$.

\section{Other problems}

This section overviews some other combinatorial problems
in which Theorem \ref{thm:main} leads to new
results.
Let $C\in H^n$ be a linear code (a linear subspace of the $\ff_2$-linear
space). Then $d(C)$ equals the minimum Hamming weight of a nonzero code
vector (the Hamming weight is the norm corresponding to the Hamming
distance $\dist(\cdot,\cdot)$). Let $\dim C=k$ and let $G$ be the
$(k\times n)$ matrix whose rows form a basis of $C$ (we assume that
$G$ has no all-zero columns). Columns of $G$ can be viewed a multiset
$X$ of points in the $(k-1)$-dimensional projective space $\pp(H^k)$; then
clearly
\[
d(C)=n-\max_{\mathcal H}|\{X\cap {\mathcal H}\}|
\]
where the maximum is taken over all hyperplanes in $\pp(H^k)$.
Likewise if $\codim {\mathcal H}=r\ge 2,$
the corresponding value is called the $r$th (higher) weight of $C$,
denoted $d_r(C).$ Properties of higher weights were a subject
of intensive study in the 1990s \cite{tsf95}. One of the problems
that present interest is finding $\max \{|C|: d_r(C)=d\}$ for a
given $r\ge 2$. Best known asymptotic upper bounds on this quantity were
proved in \cite{ash99b}, an essential ingredient of the proof
being Theorem \ref{thm:spectrum-Hamming} and related results.

Variations of the polynomial method proved to be efficient
for deriving new upper bounds on the maximum size of 
list-decodable codes \cite{ash99c}, on the covering radius
of codes with a given dual distance \cite{ash99h}, and
on the minimum distance of doubly-even self-dual codes \cite{kra99}.

The polynomial method also proved useful in the study of quantum
information transmission. It was applied in \cite{ash99e} to derive
upper bounds on the size of quantum codes. In \cite{ash99d} the concept
of error detection was extended to quantum codes.
It turned out that the probability of undetected error can be
expressed via the weight enumerators of quantum codes (the Shor-Laflamme
enumerators \cite{sho97}) in a way similar to $P_{ue}(H^n;C,p).$
It is shown in \cite{ash99d} that there exist quantum codes
for which the probability of undetected error is an exponentially
falling function of the code length $n$. Furthermore,
a version of the polynomial method developed in \cite{ash99e}, \cite{ash99d}
leads to upper bounds on this exponent which are tight in a certain
region of code rates.

\renewcommand\baselinestretch{0.9}
{\small
\providecommand{\bysame}{\leavevmode\hbox to3em{\hrulefill}\thinspace}

}
\end{document}